\documentclass[12pt]{article}

\usepackage{amsmath}
\usepackage{amssymb}

\usepackage{graphics}

\begin{document}

\title{Optimal Prediction of \\ Stiff Oscillatory Mechanics}

\author{Anton P. Kast \\
\\
Applied Numerical Algorithms Group \\
Lawrence Berkeley National Laboratory \\
Berkeley, CA  94720}

\maketitle

\begin{abstract}

We consider many-body problems in classical mechanics where a wide range of
time scales limits what can be computed.  We apply the method of optimal
prediction to obtain equations which are easier to solve numerically.  We
demonstrate by examples that optimal prediction can reduce the amount of
computation needed to obtain a solution by several orders of magnitude.

\end{abstract}

\vspace{1in}

\thanks{This work was supported in part by the Department of Energy Office
of Advanced Computing Research, Mathematical, Information, and Computational
Sciences Division, under Contract No.\ DE-AC03-76SF00098.}

\newpage

\setlength{\baselineskip}{2\baselineskip}

\section{Stiff oscillatory mechanics}

There are many problems in classical mechanics where what can be computed is
limited by the simultaneous presence of both fast and slow motion: some
variables oscillate rapidly while others change slowly, so standard numerical
methods can require a large number of time steps to give accurate answers.
Stiffness of this type limits calculations of planetary motion, drift in
high-frequency electronic oscillators, and the dynamics or large
molecules~\cite{petzold}.

For instance, in molecular dynamics it is standard~\cite{reich} to model the
motion of many atoms as a mechanical system with a Hamiltonian of the form
\begin{equation}
H = \frac{1}{2}\sum_{j=1}^N\frac{p_j^2}{2m_j}
	+ V(q_1,\ldots,q_N) + \frac{1}{2}\sum_{j=1}^N\sum_{k=1}^N
	g_j(q)A_{jk}g_k(q)
\label{cmdH}
\end{equation}
where $(q_j,p_j)$ are the coordinates and momenta of the atoms and $N$ is the
number of atoms, commonly in the range $10^4$ to $10^5$.  Here $V$ denotes a
smoothly-varying potential energy of interaction among coordinates, the $g$'s
are bond angles or interatomic spacings (functions of the coordinates), the
$m$'s are masses, and $A$ is a matrix of spring constants.  Such models are
used to describe both the large-scale motion that takes place over
milliseconds and also the rapid vibrational motions at chemical bonds which
are measured in terahertz.

In a recent paper\cite{sw}, Stuart and Warren considered a particular stiff
Hamiltonian problem of the form~(\ref{cmdH}) that was originally meant to
model a particle interacting with a heat bath~\cite{fk}, and they constructed
numerical schemes that worked well with large time steps.  They were able to
compute the motion of slowly-varying quantities accurately, even when most of
the dynamics was grossly underresolved in time (i.e., even when their time
step was much longer than the periods of most normal modes of oscillation).

This observation, that a scheme may be optimized to work well even when the
resolution is poor, is similar to the results of optimal
prediction~\cite{pnas,cpam,ams}; optimal prediction is a method for reducing
the resolution required to solve a large system of equations.  A smaller
system is constructed, designed to yield expectations of solutions of the
larger system and to be computationally practical even when the larger system
is not.  Since Stuart and Warren have found schemes for some large, stiff
systems that work with big time steps, it is natural to ask whether there are
smaller systems of differential equations (just describing the slower modes)
that would work at these big time steps.

In this paper we show how optimal prediction may be applied to a class of
large, stiff Hamiltonian systems like~(\ref{cmdH}) to yield effective
equations which are smaller and slower.  We demonstrate the method on the
Stuart-Warren model and on a generalization of it that more closely
approximates realistic models of molecular dynamics.  The benefits are longer
time steps, lower dimensionality (hence fewer force evaluations per time
step), and a systematic approach that may may be broadly applied.

\section{Optimal prediction}

Optimal prediction is a method that takes a large system of differential
equations together with a probability distribution for the dependent
variables, and produces a smaller system of equations for the expectations of
some selected variables while averaging over all the others.  The method is
described in~\cite{pnas,cpam,ams}.  Error bounds for the
method can be found in~\cite{hald}.

Suppose we are given a large dynamical system 
\begin{equation}
\dot{u}_i = R_i(u_1,\ldots,u_N), \qquad i=1,\ldots,N
\label{dynamicalSystem}
\end{equation}
for dependent variables $u_1,\ldots,u_N$, and we are also given a normalized
probability density $P(u_1,\ldots,u_N)$ which is invariant
under~(\ref{dynamicalSystem}),
\begin{equation}
\sum_{j=1}^N\frac{\partial P}{\partial u_j} R_j(u_1,\ldots,u_N) \equiv 0 .
\end{equation}

The first step in the optimal prediction procedure is to identify
``collective variables,'' meaning a small number of functions of the
dependent variables whose evolution we would like to predict.
We denote these collective variables by $v_1(u)\ldots v_n(u)$ where $n<N$.
The idea in optimal prediction is to treat the $u$'s as random, treat their
combinations in the $v$'s as known, and to estimate the rates of change of
the $v$'s by conditional expectations.

One writes out a formula for the rate of change of the $v$'s induced
by~(\ref{dynamicalSystem}),
\begin{equation}
\dot{v}_\mu(u) = \sum_{j=1}^N \frac{\partial v_\mu}{\partial u_j}
    R_j(u_1,\ldots,u_N)
\end{equation}
Then one uses $P(u)$ to compute the expectation of this expression subject to
conditions that $v_\mu(u)=\overline{v}_\mu$ for some $n$ numbers
$\overline{v}_1\ldots\overline{v}_n$,
\begin{equation}
\langle \dot{v}_\mu \rangle_{\overline{v}_1\cdots \overline{v}_n} =
\frac{
	\displaystyle{
	\int\ \dot{v}_\mu(u)\ P(u) \ 
	\prod_{\nu=1}^n \delta(v_\nu(u)-\overline{v}_\nu) \ du}}
{	\displaystyle{
	\int P(u)\ 
	\prod_{\nu=1}^n \delta(v_\nu(u)-\overline{v}_\nu) \ du}} .
\label{integrals}
\end{equation}
Finally, one hypothesizes that the mean evolution of the $v$'s is
approximated by the solutions $\overline{v}_\mu(t)$ of the new system,
\begin{equation}
\dot{\overline{v}}_\mu(t) =
	\left\langle 
\sum_{i=1}^N \frac{\partial v_\mu}{\partial u_j}
    R_j(u_1,\ldots,u_N)
	\right\rangle_{\overline{v}_1(t)\cdots\overline{v}_n(t)}
\label{hypothesis} 
\end{equation}
The new system~(\ref{hypothesis}) is a closed system of equations for the
$\overline{v}$'s, and it is $n$-dimensional instead of $N$-dimensional.

Equation~(\ref{hypothesis}) approximates the evolution of the mean values of
the $v$'s.  The idea of the approximation is that at every moment in time,
the $u$'s are distributed according to their invariant probability density
subject to conditions on the values of collective variables.  All that
changes in time is the conditions, according to our
hypothesis~(\ref{hypothesis}).  Actually, if the $v$'s were given and the
$u$'s were distributed according to a conditioned invariant distribution at
time $t=0$, then at a future time $t>0$ the $v$'s would be indeterminate and
the $u$'s would become distributed in some more general way.  Average values
of the $v$'s at all times $t>0$ would still be well-defined though, and they
would be determined by the values of the $v$'s at $t=0$.  The
system~(\ref{hypothesis}) is meant to approximate such exact mean evolutions
of collective variables from initial values.

Although equation~(\ref{hypothesis}) is conjectural, some general results are
known about its accuracy.  First, it clearly gives an asymptotically exact
prediction of mean futures for short times.  Second, it appears in an exact
formula for mean futures due to Zwanzig~\cite{zwanzig80} (recently studied by
others~\cite{chk}) which reveals corrections in terms of history integrals
and noise-like functions which are statistically uncorrelated with the
collective variables.  Third, error bounds for the method have been
established in the case of Hamiltonian dynamical systems~\cite{hald}.

There are two technical challenges in the application of~(\ref{hypothesis}):
collective variables must be selected, and the conditional expectations on
the right-hand side must be explicitly evaluated, usually requiring
approximations of the integrals in equation~(\ref{integrals}).  Both steps
are critical to accuracy.  In complex problems, therefore, the best way to
determine the usefulness of the approximation~(\ref{hypothesis}) is
empirically: one generates large random ensembles of initial conditions
for~(\ref{dynamicalSystem}), integrates each initial condition, then averages
the results to determine a mean future.  One then compares the answer to an
integral of~(\ref{hypothesis}).

In the present paper, we will consider Hamiltonian equations where the
dependent variables are canonical coordinate pairs
$(q_1,p_1)\ldots(q_N,p_N)$.  Hamiltonian equations preserve the canonical
probability density, $e^{-H}$, so we will use this as our probability
density.  We assume that the first $n$ coordinate pairs
$(q_1,p_1)\ldots(q_n,p_n)$ are of interest, and we will take the remaining
dynamical variables as random.

The optimal prediction procedure is to take the full system of Hamilton's
equations,
\begin{equation}
\dot{q}_j = \frac{\partial H}{\partial p_j}, \quad
\dot{p}_j = - \frac{\partial H}{\partial q_j}, \quad
j=1,\ldots,N,
\label{hamiltonEqns}
\end{equation}
discard the equations with indices $j>n$, and replace the right-hand sides of
the remaining equations with their expectations with respect to $e^{-H}$
conditioned by the selected variables:
\begin{equation}
\dot{q}_\mu = \left\langle \frac{\partial H}{\partial p_\mu}\right\rangle_n, \quad
\dot{p}_\mu = \left\langle - \frac{\partial H}{\partial q_\mu}\right\rangle_n, \quad
\mu=1,\ldots,n
\label{opEqns}
\end{equation}
where $\langle\cdot\rangle_n$ denotes the conditioned expectation,
\begin{equation}
\langle f \rangle_n = Z^{-1} \int \prod_{j=n+1}^N dq_j \ dp_j \ e^{-H}
	f(q_1,\ldots,q_N;p_1,\ldots,p_N)
\label{condDiff}
\end{equation}
with $Z$ a normalization constant.  For any function $f$ of the canonical
variables, $\langle f \rangle_n$ is a function of $q_1\cdots q_n$, $p_1\cdots
p_n$ only, so the $2n$-dimensional system of equations~(\ref{opEqns}) is
closed.

The reduced system~(\ref{opEqns}), the first approximation in optimal
prediction, defines an approximate solution to a Liouville problem for the
evolution of a probability measure on phase space.  At least for short times,
the system~(\ref{opEqns}) is guaranteed to give the expectations of the
selected variables, averaging over all possible initial data for the
discarded variables.

We need to evaluate the conditional expectations in~(\ref{opEqns}).  This is
easy if $e^{-H}$ is a Gaussian distribution (i.e., if $H$ is quadratic, or
equivalently if the equations of motion are linear).  If $e^{-H}$ is not
Gaussian, perturbative techniques are available to approximate its
expectations by Gaussian expectations.  Thus the following results for
Gaussian distributions will be sufficient for our purposes,
see~\cite{pnas,cpam,ams} for details.

Let $x_1,\ldots,x_N$ be Gaussian random variables distributed with density
\begin{equation}
P(x_1,\ldots,x_N) \propto \exp\left(-\frac{1}{2}\sum_{j=1}^N\sum_{k=1}^N
	x_jA_{jk}x_k + \sum_{j=1}^N b_j x_j \right) .
\end{equation}
We denote expectations with respect to this density by $\langle\cdot\rangle$,
and $\langle x_i \rangle = \sum_{j=1}^N A^{-1}_{ij}b_j$.  Now suppose that
$x_1 \ldots x_n$ are given for all $n < N$.  The conditional expectations
of $x_{n+1} \ldots x_N$ conditioned by $x_1 \ldots x_n$ are denoted $\langle
x_i \rangle_n$, $i=n+1,\ldots,N$ and are given explicitly by
\begin{equation}
\langle x_i \rangle_n = \langle x_i \rangle +
\sum_{\mu=1}^n \sum_{\nu=1}^n A^{-1}_{i\mu} M^{-1}_{\mu\nu}
(x_\nu - \langle x_\nu \rangle), \qquad i=n+1,\ldots,N
\label{constrainedAv}
\end{equation}
where $M_{\mu\nu}=A^{-1}_{\mu\nu}$ for $\mu,\nu=1,\ldots,n$ and $M^{-1}$ is
the inverse of the $n\times n$ (not $N \times N$) matrix $M$.

The conditioned covariances, $\text{Cov}_n(x_i,x_j)=\langle
x_ix_j\rangle_n-\langle x_i\rangle_n\langle x_j\rangle_n$ are given in terms
of the unconditioned expectations $\text{Cov}(x_i,x_j)=\langle
x_ix_j\rangle-\langle x_i\rangle\langle x_j\rangle$ by
\begin{equation}
\text{Cov}_n(x_i,x_j) = \text{Cov}(x_i,x_j)
	- \sum_{\mu=1}^n \sum_{\nu=1}^n
	A^{-1}_{i\mu} M^{-1}_{\mu\nu} A^{-1}_{\nu j} .
\label{constrainedCov}
\end{equation}

The conditioned expectation of any polynomial in $x_1\ldots x_N$ may be found
from these formulae by Wick's theorem.

\section{Generalizations of the Stuart-Warren experiments}

Stuart and Warren~\cite{sw} (see also~\cite{cssw},~\cite{fk},
and~\cite{ford}) considered a one-dimensional collection of particles
connected by springs.  There was one distinguished particle with mass $1$,
coordinate $Q$ and momentum $P$.  The distinguished particle was connected by
springs of spring constant $k$ to $N$ other particles with masses $k/j^2$,
coordinates $q_j$ and momenta $p_j$, $j=1\ldots N$, representing a heat bath.

The motion of this collection of particles and springs is defined by the
Hamiltonian
\begin{multline}
H(Q,P;q_1,\ldots,q_N;p_1,\ldots,p_N) \\
= \frac{1}{2}(V(Q)+P^2)
+ \sum_{j=1}^N \left[ \frac{p_j^2}{2m_j} + \frac{1}{2}k(Q-q_j)^2 \right]
\label{originalH}
\end{multline}
where $(Q,P)$ and $(q_j,p_j)$ are canonically conjugate dynamical variables
for $j=1,\ldots,N$ and $m_j=k/j^2$.  The equations of motion are 
\begin{equation}
\begin{aligned}[2]
\dot{Q} &= P & \dot{P} &= - V'(Q) + k \sum_{j=1}^N (q_j-Q) \\
\dot{q}_j &= p_j/m_j \quad & \dot{p}_j &= k ( Q - q_j ), \quad j=1,\ldots,N
\end{aligned}
\label{original}
\end{equation}

This system is of the form~(\ref{cmdH}) (with an extra pair of coordinates
$(Q,P)$), and it is chosen so that fast and slow motion are separated:
lighter particles will move faster, heavier particles will move slower, and
the mass $m_j$ goes down as $j$ goes up.

A central result of~\cite{sw} is that if all the heat bath particles start
out randomly, with statistics determined by the canonical distribution, then
in the limit $N\rightarrow\infty$ the coordinate of the distinguished
particle obeys the stochastic equation,
\begin{equation}
\ddot{Q} + \frac{k\pi}{2}\dot{Q} + V'(Q) - \frac{k}{2} Q = F
\end{equation}
where $F(t)$ is a stochastic process related to white noise.  This equation
for $Q$ is remarkable because it makes no reference to the history of
$Q$---it is a differential equation, not an integro-differential equation.
In a general Hamiltonian problem, if one variable $Q$ is fixed initially and
the others are random, at future times there is no time-invariant
relationship among the expectation of $Q$ and its time
derivatives~\cite{mori1,mori2,zwanzig}.  The first approximation of optimal
prediction~(\ref{opEqns}) may be characterized as the assumption that the
values of the selected variables do determine their own future expectations.
In general this assumption is not exactly true, but in the Stuart-Warren
model it is true exactly in the $N\rightarrow\infty$ limit.

Stuart and Warren proceeded to integrate their model with large time steps.
If $Q$ were fixed, then each $q_j$ would oscillate harmonically with
frequency $\omega_j=j$.  This implies that a discretization of the $2N+2$
equations~(\ref{original}) would be resolved in time if $N\Delta t \ll 1$.
If this condition on $\Delta t$ were violated, then the result of the
computation would depend on how the equations were discretized.  The
intriguing result of~\cite{sw} is that some schemes will give the right
evolution for $Q$ and $P$ when $N\Delta t \gtrsim 1$ and others will not.
For instance, if the scheme is
\begin{equation}
\begin{aligned}[2]
\frac{Q^{n+1}-Q^n}{\Delta t} &= P^{n+1} &
\frac{P^{n+1}-P^n}{\Delta t} &= - V(Q^n) + k \sum_{j=1}^N (q_j^{n+\sigma}-Q^n) \\
\frac{q_j^{n+1}-q_j^n}{\Delta t} &= p_j^{n+1} / m_j \quad &
\frac{p_j^{n+1}-p_j^n}{\Delta t} &= k(Q^n-q_j^n) \quad j=1,\ldots,N
\end{aligned}
\label{swScheme}
\end{equation}
then $\sigma=0$ (a symplectic method) gives the right answer for $Q$ and $P$,
but $\sigma=1$ (another convergent method) does not.

For concreteness, we pick $V(Q)=\frac{1}{2}Q^2$.  Since $H$
in~(\ref{originalH}) is then quadratic, the canonical probability density is
Gaussian, and formula~(\ref{constrainedAv}) gives the conditioned
expectations as
\begin{equation}
\langle q_j \rangle_n = Q, \quad \langle p_j \rangle_n = 0 \quad (n<j\leq N) .
\end{equation}
Taking the conditional expectations of the right-hand sides
of~(\ref{original}) and evaluating them using these results, we find that
the equations of optimal prediction are
\begin{equation}
\begin{aligned}[2]
\dot{Q} &= P & \dot{P} &= - Q + k \sum_{\mu=1}^n (q_\mu-Q) \\
\dot{q}_\mu &= p_\mu/m_\mu \quad & \dot{p}_\mu &= k ( Q - q_\mu ), \quad \mu=1,\ldots,n
\end{aligned}
\label{originalOP}
\end{equation}
These are identical in form to the original equations~(\ref{original}).  It
comes as no surprise, therefore, that the motion of $Q$ can be computed with
large $\Delta t$: pick the $\Delta t$ desired, find an $n\ll N$ such that
$n\Delta t \ll 1$, and perform a resolved integration of~(\ref{originalOP})
with this $n$ and $\Delta t$.  Reasonable approximations for the selected
variables are guaranteed, at least for short times.

\begin{figure}
\resizebox{\textwidth}{!}
	{\includegraphics*{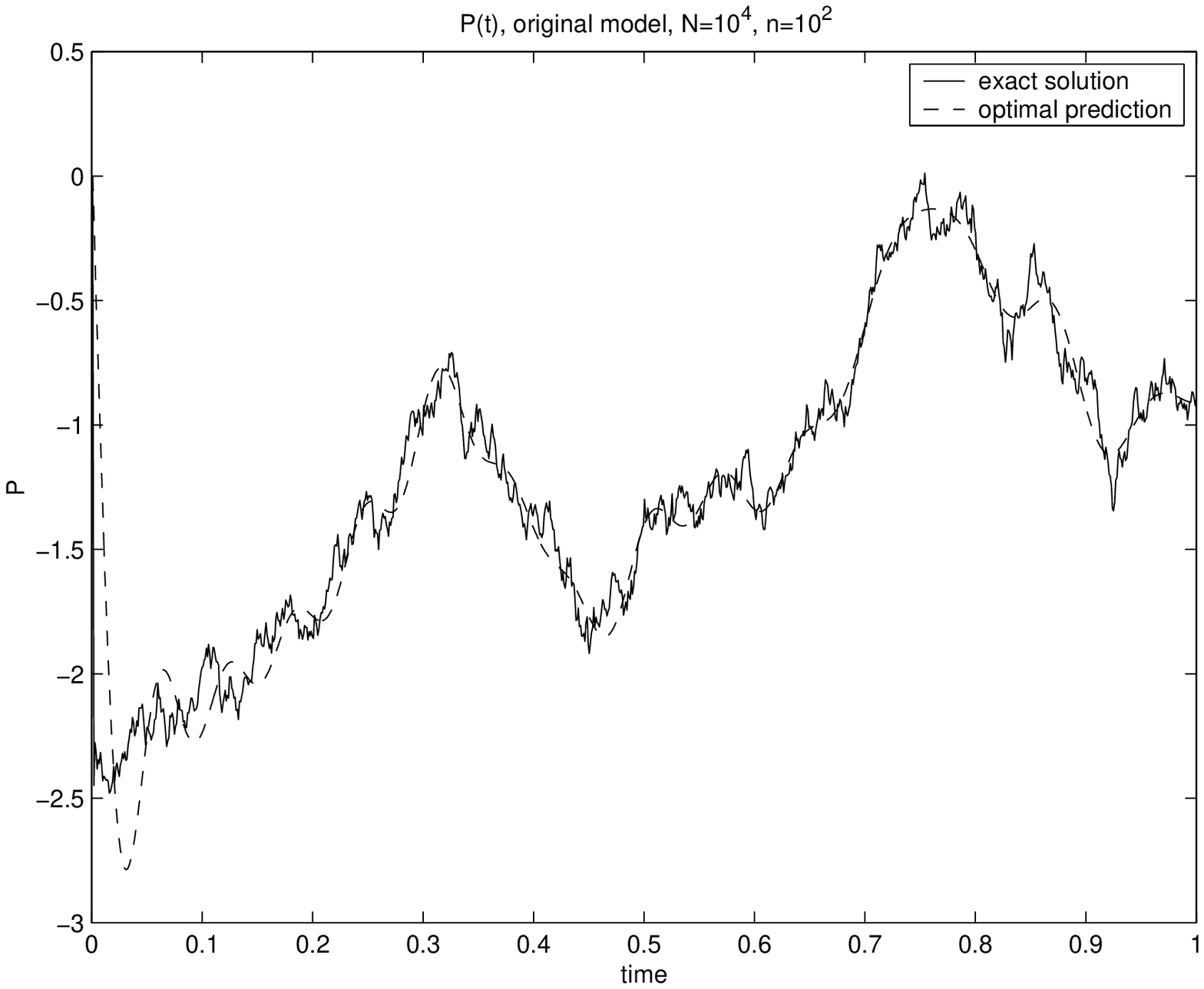}}
\caption{The evolution of $P(t)$ determined in two ways: by solving the
equations of motion~(\ref{original}) with $N=10,000$ particles and random
initial data (exact evolution, $\Delta t=10^{-2}/N$); and by solving the
reduced equations~(\ref{reduced}) with $n=100$ particles and a time step $100$
times longer (optimal prediction, $\Delta t=1/N=10^{-2}/n$).  For these
calculations, $k_Q=k_q=1$.}
\label{f1}
\end{figure}

Figure~\ref{f1} shows a fully-resolved calculation ($N\Delta t = 10^{-2}$) of
$P(t)$ starting from $P(0)=0$, $Q(0)=1.5$, with $q_j(0)$ and $p_j(0)$ chosen
randomly from the canonical ensemble (i.e., chosen with probability density
$e^{-H}$) conditioned by $Q(0)$ and $P(0)$.  It also shows the solution to
the same problem as computed by a resolved integration of~(\ref{originalOP}),
which was achieved with $n\Delta t = 10^{-2}$.  The optimal prediction
calculation accurately duplicates the low-frequency behavior of the exact
solution, and it does so in fewer dimensions with a larger time step.  In
this case, with $N=10^4$ and $n=10^2$, the optimal prediction curve was
about $10,000$ times faster to compute than the resolved solution.  The
optimal prediction has the further advantage that it did not use the initial
data $q_{n+1}(0) \ldots q_N(0)$, $p_{n+1}(0) \ldots p_N(0)$ and may claim to
be an average answer over all possible values of these data.

\section{More general models}

Realistic applications, such as molecular dynamics, involve more complex
interactions than are present in the model~(\ref{original}).  In particular,
we may expect that every particle would interact with every other, and that
the interactions would be nonlinear.

We therefore consider a generalization of the model~(\ref{original}) where
every $q_1\ldots q_N$ is coupled to every other $q_1\ldots q_N$ by a spring,
and the springs are nonlinear:  
\begin{multline}
H(q_1,\ldots,q_N;p_1,\ldots,p_N) \\
= \sum_{j=1}^N \frac{p_j^2}{2m_j}
+ \frac{1}{2} k^{(2)} \sum_{j=1}^N \sum_{l=j+1}^N (q_j-q_l)^2 
+ \frac{1}{4} k^{(4)} \sum_{j=1}^N \sum_{l=j+1}^N (q_j-q_l)^4
\end{multline}
\begin{equation}
\left.
\begin{aligned}
\dot{q}_j &= p_j/m_j \\
\dot{p}_j &= - k^{(2)} \sum_{l=1}^N (q_j-q_l)
			 - k^{(4)} \sum_{l=1}^N (q_j-q_l)^3
\end{aligned}
\right\} j=1,\ldots,N .
\label{new}
\end{equation}
This model makes no reference to a distinguished particle; each one of the
$N$ particles interacts with all of the others through the same potential
energy, which is parameterized by the new spring constants $k^{(2)}$ and
$k^{(4)}$.

We derive the optimal prediction equations of the system~(\ref{new}) for $q_1
\ldots q_n$, $p_1 \ldots p_n$ by averaging over $q_{n+1} \ldots q_N$,
$p_{n+1} \ldots p_N$.  Since the interactions are now nonlinear, the
probability density $e^{-H}$ is no longer Gaussian, so we must work harder to
evaluate the conditioned expectations.

Hald has observed, as reported in~\cite{chk}, that optimal prediction
equations of the form~(\ref{opEqns}) are always Hamiltonian, and that their
Hamiltonian is
\begin{equation}
H'(q_1,\ldots,q_n;p_1,\ldots,p_n) = -\log\left(
\int\prod_{j=n+1}^N\,dq_j\,dp_j\,e^{-H} \right).
\label{effH}
\end{equation}
We may therefore approximate the conditioned expectations of~(\ref{opEqns})
by first approximating $H'$, and then deriving~(\ref{opEqns}) by
differentiation:
\begin{equation}
\dot{q}_\mu=\frac{\partial H'}{\partial p_\mu}, \qquad
\dot{p}_\mu=-\frac{\partial H'}{\partial q_\mu}, \qquad \mu=1,\ldots,n .
\label{diffEffH}
\end{equation}
We decompose $H$ into its quadratic part plus its higher-order part,
\begin{equation}
\begin{aligned}
H&=H_0+H_1 \\
H_0&=\sum_{j=1}^N\frac{p_j}{2m_j} +
	\frac{k^{(2)}}{2}\sum_{j=1}^N\sum_{l=j+1}^N(q_j-q_l)^2 \\
H_1&=\frac{k^{(4)}}{4}\sum_{j=1}^N\sum_{l=j+1}^N(q_j-q_l)^4
\end{aligned}
\end{equation}
and proceed by determining $H'$ perturbatively as a power series in
$k^{(4)}$.  An alternate method for perturbative treatment of optimal
prediction is described in~\cite{clk}.

Hald's formula~(\ref{effH}) implies
\begin{equation}
\begin{aligned}
H'&=-\log\left(\int\prod_{j=n+1}^N\,dq_j\,dp_j\,e^{-H_0}\right)
   - \log\left(\frac{\int\prod_{j=n+1}^n\,dq_j\,dp_j\,e^{-H_0}e^{-H_1}}
                    {\int\prod_{j=n+1}^N\,dq_j\,dp_j\,e^{-H_0}}\right) \\
  &= \text{($H_0$-part)} - \log\left\langle e^{-H_1}\right\rangle_{n,0}
\end{aligned}
\label{haldImplies}
\end{equation}
where the new average, $\langle\cdot\rangle_{n,0}$ denotes an average with
respect to the conditioned {\em Gaussian} measure, defined just as in the
definition~(\ref{condDiff}) but with $H_0$ replacing $H$.  The
``($H_0$-part)'' term would be the effective Hamiltonian if $H_1$ were zero,
and it contributes linear terms to the equations of motion which are easily
evaluated by the regression formula~(\ref{constrainedAv}).  The other term
in~(\ref{haldImplies}) is equal to a power series in $k^{(4)}$,
\begin{equation}
\log\langle e^{-H_1}\rangle_{n,0} = \sum_{m=1}^\infty \frac{(-1)^m}{m!}
\langle H_1^m\rangle_{n,0}^{(c)}
\end{equation}
where $\langle H_1^m\rangle_{n,0}^{(c)}$ denotes the $m$-th cumulant of $H_1$ with
respect to the conditioned Gaussian measure.  Each cumulant in this series
may be evaluated by Wick's theorem, where only ``connected'' pairings (in the
sense of perturbation theory in physics) are included.

To first order in $k^{(4)}$, we need to evaluate
\begin{equation}
\begin{aligned}
\langle H_1^1 \rangle_{n,0}^{(c)} &= \langle H_1 \rangle_{n,0} \\
&= \left\langle \frac{k^{(4)}}{4} \sum_{j=1}^N\sum_{l=j+1}^N
   (q_j-q_l)^4\right\rangle_{n,0} \\
&= \frac{k^{(4)}}{4}\left[
   \sum_{\mu=1}^n\sum_{\nu=\mu+1}^n(q_\mu-q_\nu)^4
   +
   \sum_{\mu=1}^n\sum_{l=\mu+1}^n\left\langle(q_\mu-q_l)^4\right\rangle_{n,0} \right]
   + \text{(constant)}
\end{aligned}
\end{equation}
where ``$\text{(constant)}$'' denotes terms that are independent of
$q_1,\ldots,q_n$ and $p_1,\ldots,p_n$ (and therefore do not affect equations
of motion).  The average $\langle\cdot\rangle_{n,0}$ may be deduced from the
expectations,
\begin{equation}
\begin{aligned}
\langle q_j \rangle_{n,0} &= \frac{1}{n} \sum_{\mu=1}^n q_\mu \\
\text{Cov}_0(q_j,q_l) &=\frac{1}{Nk^{(2)}}(1+\delta_{jl})
\end{aligned}
\qquad j,l=n+1,\ldots,N
\end{equation}
together with Wick's theorem.  The result for $H'$, to first order in
$k^{(4)}$, is
\begin{equation}
\begin{aligned}
H' = \sum_{\mu=1}^n\frac{p^2_\mu}{2m_\mu}
   &+ \frac{C_2}{2}\sum_{\mu=1}^n\sum_{\mu=\nu+1}^n(q_\mu-q_\nu)^2 \\
   &+ \frac{C_4}{4}\sum_{\mu=1}^n\sum_{\mu=\nu+1}^n(q_\mu-q_\nu)^4 \\
   &+ \frac{D_4}{4}\sum_{\mu=1}^n\left(q_\mu-\frac{1}{n}\sum_{\nu=1}^nq_\nu\right)^4
   + O\left(k^{(4)}\right)^2
\label{answer}
\end{aligned}
\end{equation}
where the coupling constants to this order in $k^{(4)}$ are
\begin{equation}
\begin{aligned}
C_2 &= \frac{N}{n}k^{(2)} + 3 \frac{(N-n)(n+1)}{Nn} \frac{k^{(4)}}{k^{(2)}} \\
C_4 &= k^{(4)} \\
D_4 &= k^{(4)} (N-n)
\end{aligned} .
\end{equation}
We differentiate~(\ref{answer}) to obtain the optimal prediction equations
for the new system~(\ref{new}) to $O(k^{(4)})^2$,
\begin{equation}
\begin{aligned}
\dot{q}_\mu &= p_\mu/m_\mu \\
\dot{p}_\mu &= - C_2 \sum_{\nu=1}^n(q_\mu-q_\nu)
               - C_4 \sum_{\nu=1}^n(q_\mu-q_\nu)^3 \\
			&\phantom{=} - D_4 \frac{1}{n}\sum_{\nu=1}^n
			   \left[\left(q_\mu-\frac{1}{n}\sum_{\sigma=1}^nq_\sigma\right)^3
			       - \left(q_\nu-\frac{1}{n}\sum_{\sigma=1}^nq_\sigma\right)^3
			   \right]
\end{aligned}
\qquad
\mu=1,\ldots,n .
\label{reduced}
\end{equation}

\begin{figure}
\resizebox{\textwidth}{!}
	{\includegraphics*{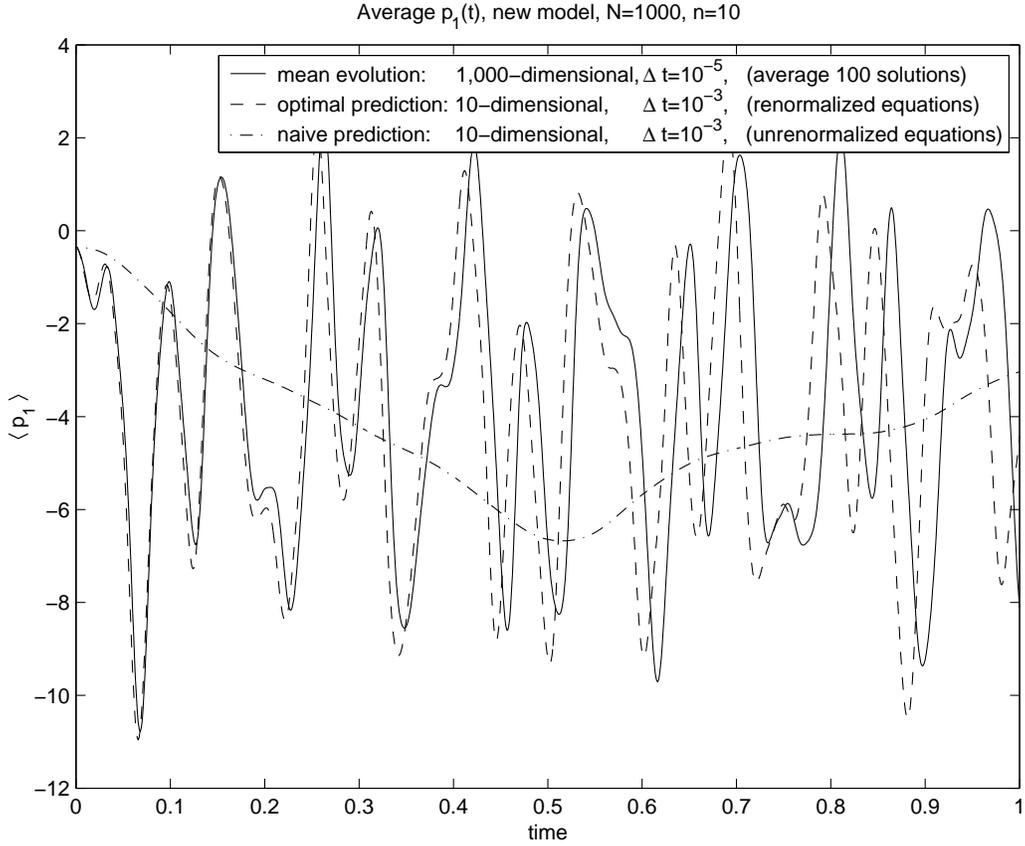}}
\caption{The {\em average} evolution of $p_1(t)$ determined in three ways: by
solving the equations of motion~(\ref{new}) for $100$ different initial
conditions, with $N=10^3$ particles, $\Delta t=10^{-2}/N$, and then averaging
all $100$ solutions (mean evolution); by solving the reduced
equations~(\ref{reduced}) once, with $n=10$ particles, $\Delta t = 1/N =
10^{-2}/n$ (optimal prediction); and by solving the original
equations~(\ref{new}) once with $N=10$, $\Delta t=10^{-2}/N$ (naive
prediction, just neglecting interactions with discarded variables).}
\label{f2}
\end{figure}

We performed a more rigorous test of this new model, comparing it to an
actual mean evolution.  The results are shown in Figure~\ref{f2}.  We once
again picked $q_1 \ldots q_n$, $p_1 \ldots p_n$ ($n=10$) from the canonical
distribution $e^{-H}$ for $N$ particles ($N=1000$ at $k^{(2)}=1$ and
$k^{(4)}=0.1$).  We then generated an ensemble of $100$ sets of values for
$q_{n+1} \ldots q_N$, $p_{n+1} \ldots p_N$ from the canonical distribution
conditioned by $q_1\ldots q_n$, $p_1\ldots p_n$, and for each set integrated
the equations~(\ref{new}).  Averaging over all $100$ solutions yielded the
solid curve for $p_1(t)$.  We then discarded the ensemble and used the
original $q_1 \ldots q_n$, $p_1 \ldots p_n$ as initial conditions for the
reduced system~(\ref{reduced}), which we integrated with $\Delta
t=10^{-2}/n=1/N$.  This $\Delta t$ is small enough to resolve the reduced
dynamics but much too large to resolve the original dynamics.  The solution
for $p_1(t)$ from~(\ref{reduced}) is the dashed curve.  Finally, for
comparison we performed the naive experiment of simply truncating the big
system~(\ref{new}) to $n$ degrees of freedom, effectively ignoring the
lighter particles without changing the interactions.  This produced the
dot-dashed curve.

The figure shows that the reduced system accurately predicts the average
evolution of $p_1(t)$, and it does so with $1$ percent of the degrees of
freedom and time steps that are $100$ times larger.  The naive experiment
shows that the new couplings are critical to the answer.  Since forces must
be evaluated $N(N-1)/2$ times per time step for $N$ particles, optimal
prediction speeds up the calculation of $p_1(t)$ in this case by about a
factor of about $10^6$.

\section{Conclusions}

We have shown that optimal prediction may be applied to large, stiff
Hamiltonian systems of differential equations to make new systems that
are smaller, better-conditioned, and approximate the original equations in
the mean.  We have demonstrated that the method gives accurate answers while
allowing larger time steps and requiring fewer force evaluations.

\section{Acknowledgements}

The author thanks Profs. A.\ Chorin, O. Hald, R. Kupferman, and A.\ Stuart
for helpful discussions.

\end{document}